\documentclass[hidelinks]{siamart220329}

\usepackage{amsmath,amsfonts,amssymb}
\usepackage{latexsym}
\usepackage{graphicx}
\usepackage{tikz}
\usepackage{hyperref}
\usepackage{colonequals}

\newcommand{\RR}{\mathbb R}

\DeclareGraphicsExtensions{.pdf,.eps,.ps,.eps.gz,.ps.gz,.eps.Y}

\title{Degenerate area preserving surface Allen-Cahn equation and its sharp interface limit}

\author{Michal Bene{\v s}\thanks{Department of Mathematics, Faculty of Nuclear Sciences and Physical Engineering, Czech Technical University in Prague, 12000 Prague, Czech Republic (\email{michal.benes@fjfi.cvut.cz}, \email{miroslav.kolar@fjfi.cvut.cz})}
\and Miroslav Kol{\'a}{\v r}\footnotemark[1]
\and Jan Magnus Sischka\thanks{Institut f\"ur Wissenschaftliches Rechnen, Technische Universit\"at Dresden, 01062 Dresden, Germany (\email{jan\_magnus.sischka@tu-dresden.de}, \email{axel.voigt@tu-dresden.de})
\funding{JMS and AV were supported by the German Research Foundation
(DFG) through EXC PoL. MB and MK were partly supported by the project 
{\it 21-09093S} of the Czech Science Foundation. We further acknowledge computing resources provided by ZIH
at TU Dresden within projects WIR.}}
\and Axel Voigt\footnotemark[2]}

\headers{Degenerate area preserving surface Allen-Cahn equation}{Michal Bene{\v s}, Miroslav Kol{\'a}{\v r}, Jan M. Sischka, and Axel Voigt}

\begin{document}

\maketitle
\begin{abstract}
We consider formal matched asymptotics to show the convergence of a degenerate area preserving surface Allen-Cahn equation to its sharp interface limit of area preserving geodesic curvature flow. The degeneracy results from a surface de Gennes-Cahn-Hilliard energy and turns out to be essential to numerically resolve the dependency of the solution on geometric properties of the surface. We experimentally demonstrate convergence of the numerical algorithm, which considers a graph formulation, adaptive finite elements and a semi-implicit discretization in time, and uses numerical solutions of the sharp interface limit, also considered in a graph formulation, as benchmark solutions.
\end{abstract}



\section{Introduction}

The connection between phase field approximations and geometric partial differential equations is well established and can formally be justified by matched asymptotics, see \cite{Fifeetal_EJDE_1995}. Geometric partial differential equations are evolution equations that evolve curves or surfaces according to their curvature. Similarly to these curvature driven flows in 2D or 3D one can consider the evolution of curves on surfaces. The evolution of these curves is governed by geodesic curvature and thus in addition strongly depending on the local geometric properties of the underlying surface. First analytical attempts to connect these geodesic evolution laws to surface phase field models have been considered in  \cite{Elliottetal_SIAMJAM_2010,GKRR_MMMAS_2016,R_AML_2016,EHS_IFB_2022}. We here show this connection for a degenerate area preserving surface Allen-Cahn equation and an area preserving geodesic curvature flow. The surface Allen-Cahn equation provides the mathematical basis to study curvotaxis of cells in epithelia tissue \cite{HWV_EPL_2022}.

We consider a surface de Gennes-Cahn-Hilliard energy
\begin{equation}
{\cal{F}}_{dGCH}(\phi) = \tilde{\sigma} \int_{\cal{S}} \frac{1}{G(\phi)} \left(\frac{\epsilon}{2} \|\nabla_{\cal{S}} \phi\|^2 + \frac{1}{\epsilon} W(\phi) \right) d{\cal{S}}
\label{eq:FdGCH}
\end{equation}
with surface ${\cal{S}}$, phase field variable $\phi$, surface gradient $\nabla_{\cal{S}}$, double well potential $W(\phi) = \frac{1}{4}(\phi^2 - 1)^2$, rescaled surface tension $\tilde{\sigma}$ and small parameter $\epsilon > 0$ determining the thickness of the diffuse interface. The factor $1/G(\phi)$ is called the de Gennes coefficient in polymer science. We consider $G(\phi) = \frac{3}{2} (\phi^2 -1)$ or a regularized version $G_\eta(\phi) = (\frac{9}{4} (\phi^2 -1)^2 + \eta^2 \epsilon^2)^{1/2}$ with $\eta > 0$. The scaling coefficient is such that the sharp interface limit equals the one obtained from the usual Cahn-Hilliard energy without the de Gennes coefficient \cite{SVW_MMAS_2021,DRW_arXiv_2022}. Evolution equations based on this energy, at least in flat space, have been shown numerically advantageous, as the singularity, $G(\phi)$ or $G_\eta(\phi)$, helps to keep solutions confined in $[-1,1]$. Even though the theoretical foundation of this argument remains open, several numerical studies confirm this \cite{Naffouti2017,AEVT_PRM_2020}. The resulting degenerate area preserving surface Allen-Cahn equation on ${\cal{S}}$ reads
\begin{eqnarray}
    \label{eq:AC}
    \epsilon \tilde{\beta} G(\phi) \partial_t \phi &=& \tilde{\sigma} \left( \epsilon \Delta_{\cal{S}} \phi - \frac{1}{\epsilon} W^\prime(\phi)\right) + G(\phi) \lambda, \\
    \label{eq:AC_lambda}
    \lambda &=& \tilde{\sigma} \frac{1}{|{\cal{S}}|} \int_{\cal{S}} \left( - \frac{\epsilon}{G(\phi)}\Delta_{\cal{S}} \phi + \frac{1}{\epsilon G(\phi)} W^\prime(\phi) \right) \, d {\cal{S}},
\end{eqnarray}
with initial condition $\phi(0) = \phi_0$. $\Delta_{\cal{S}}$ denotes the Laplace-Beltrami operator and $\tilde{\beta} > 0$ is a rescaled kinetic coefficient. $\lambda$ is a Lagrange multiplier for the area constraint $\frac{1}{|{\cal{S}}|} \int_{\cal{S}} \phi \; d {\cal{S}} = \alpha$ with $\alpha \in [-1,1]$. The model follows from eq. \cref{eq:FdGCH} as a constrained $L^2$-gradient flow and considers the asymptotic approximation $\frac{\epsilon}{2} |\nabla_{\cal{S}} \phi|^2 \approx \frac{1}{\epsilon} W(\phi)$, see \cite{SVW_MMAS_2021}. The idea for this approximation  was first used for phase field approximations of surface diffusion in \cite{RRV_JCP_2006}, where $G(\phi)$ was introduced as a stabilizing function. The zero-levelset of the solution of eqs. \eqref{eq:AC} and \eqref{eq:AC_lambda} provides an approximation of a curve $\gamma$ evolving on the surface ${\cal{S}}$.

The area preserving geodesic curvature flow for the curve $\gamma$ on ${\cal{S}}$ reads
\begin{equation}
    \beta {\cal{V}} = - \sigma {\cal{H}}_\gamma + \sigma \frac{1}{|\gamma|} \int_\gamma {\cal{H}}_\gamma \; d \gamma
    \label{eq:MCF}
\end{equation}
with initial condition $\gamma(0) = \gamma_0$. ${\cal{V}}$ is the velocity of $\gamma(t)$ in the direction of the co-normal $\mu$ and ${\cal{H}}_\gamma$ is the geodesic curvature of $\gamma$. We consider ${\cal{S}}= {\cal{S}}^1(t) \cup \gamma(t) \cup {\cal{S}}^2(t)$ and enforce the constraint $|{\cal{S}}^1(t)| - |{\cal{S}}^2(t)| + \alpha |{\cal{S}}| = 0$ for each $t \in [0,T]$. The kinetic coefficient and line tension are related by $\beta = \frac{4 \sqrt{2}}{5} \tilde{\beta}$ and $\sigma = \frac{2 \sqrt{2}}{3} \tilde{\sigma}$.

Besides the connection between eqs. \eqref{eq:AC} and \eqref{eq:AC_lambda} and eq. \eqref{eq:MCF} by formal matched asymptotics, we use numerical solutions of eq. \eqref{eq:MCF} in a graph formulation, see \cite{Kolaretal_DCDSB_2017}, as benchmark problems for a numerical approach to eqs. \eqref{eq:AC} and \eqref{eq:AC_lambda}, again using a graph formulation. We use adaptive finite elements to discretize in space and a semi-implicit time-stepping scheme.

\section{Matched asymptotic analysis}

We closely follow \cite{EHS_IFB_2022} in the analysis of a phase field model in the context of two-phase biomembranes, and use the tools introduced in \cite{Elliottetal_SIAMJAM_2010} to extend the formal matched asymptotics for the area preserving Allen-Cahn equations in flat space \cite{Rubensteinetal_IMAJAM_1992} and for the de Gennes-Cahn-Hilliard energy in flat space \cite{SVW_MMAS_2021} to surfaces. We demonstrate that eqs. \eqref{eq:AC} and \eqref{eq:AC_lambda} formally converges to eq. \eqref{eq:MCF} for $\epsilon \to 0$. We therefore require $\gamma(t)$ to be a $C^1$ closed curve.

\subsection{Expansions and matching conditions}
By $(\phi_\epsilon,\lambda_\epsilon)$, we denote a family of solutions of eqs. \eqref{eq:AC} and \eqref{eq:AC_lambda} that converge formally to some limit denoted by $(\phi,\lambda)$. We assume that $\phi = \chi_\gamma$ with $\chi_\gamma : BV({\cal{S}}) \to \RR$ with $\chi_\gamma = -1$ on ${\cal{S}}^1$ and $\chi_\gamma = 1$ on ${\cal{S}}^2$ for some smooth curve $\gamma$ that separates the regions ${\cal{S}}^1 = \{(x,t) \in {\cal{S}} \times [0,T] : \phi(x,t) = -1\}$ and ${\cal{S}}^2 = \{(x,t) \in {\cal{S}} \times [0,t] : \phi(x,t) = 1\}$. We consider an outer and an inner expansion
\begin{eqnarray}
    f_\epsilon(x,t) &=& f_0(x,t) + \epsilon f_1(x,t) + \epsilon^2 f_2(x,t) + \ldots \\
    F_\epsilon(x,t) &=& F_0(s,z,t) + \epsilon F_1(s,z,t) + \epsilon^2 F_2(s,z,t) + \ldots,
\end{eqnarray}
respectively, with $f_\epsilon(x,t) = F_\epsilon(s,z,t)$, $z = r/\epsilon$ and $\Theta(s,r,t)$ a parametrization such that $s \to \Theta(s,0,t)$ is a paramtetrization of $\gamma(t)$ on ${\cal{S}}$ and $r$ is the signed geodesic distance of $x = \Theta(s,r,t) \in {\cal{S}}$ to $\gamma(t)$. Thereby, $f_\epsilon = \phi_\epsilon, \lambda_\epsilon$ and $F_\epsilon = \Phi_\epsilon, \Lambda_\epsilon$ and the outer expansion holds away from $\gamma(t)$ and the inner expansion near $\gamma(t)$. In regions where both expansions are valid the matching conditions hold
\begin{eqnarray}
    F_0(s, \pm \infty, t) &=& f_0^\pm(x,t), \label{eq:m1}
    \\ \partial_z F_0(s, \pm \infty, t) &=& 0, \label{eq:m2}
    \\ \partial_z F_1(s, \pm \infty, t) &=& \nabla_{\cal{S}} f_0^\pm(x,t) \cdot \mu(x,t), \label{eq:m3}
\end{eqnarray}
with $f_0^\pm(x,t) = \lim_{\delta\to 0} f(\Theta(s,\pm\delta,t),t)$.

\subsection{Outer solution}

Considering the terms of $O(\epsilon^{-1})$ in eq. \eqref{eq:AC} leads to
$W^\prime(\phi_0) = 0$ and thus
\begin{equation} \label{eq:o_-1}
    \phi_0 = \pm 1.
\end{equation}

\subsection{Inner solution}

The Laplace-Beltrami operator and the time derivative of $F_\epsilon(s,z,t)$ gives
\begin{eqnarray}
\Delta_{\cal{S}} F_\epsilon &=& \frac{1}{\epsilon^2} \partial_{zz} F_\epsilon + \frac{{\cal{H}}_\gamma}{\epsilon} \partial_z F_\epsilon + \partial_{ss} F_\epsilon + O(\epsilon) \label{eq:delta} \\
\frac{d}{d t}F_\epsilon &=& - \frac{1}{\epsilon} {\cal{V}} \partial_z F_\epsilon + \partial_t F_\epsilon  + O(\epsilon) \label{eq:time}
\end{eqnarray}
with $\partial_{ss}$ the second derivative along $\gamma$. Considering terms of $O(\epsilon^{-1})$ in eq. \eqref{eq:AC} leads to
$0 = \tilde{\sigma} \left( \partial_{zz} \Phi_0 - W^\prime(\Phi_0) \right)$. Using the outer expansion eq. \eqref{eq:o_-1} and the matching condition eq. \eqref{eq:m1} shows that $\Phi_0(z,s,t)$ is a solution of $\partial_{zz} \Phi_0 = W^\prime(\Phi_0)$ with $\Phi_0(\pm \infty) = \pm 1$ and thus
\begin{equation} \label{eq:tanh}
    \Phi_0(z) = \tanh \left(\frac{z}{\sqrt{2}} \right)
\end{equation}
independent of $s$ and $t$. Using this in $O(\epsilon^0)$ of eq. \eqref{eq:AC} leads to
\begin{equation}
    - \tilde{\beta} {\cal{V}} G(\Phi_0) \partial_z \Phi_0 = \tilde{\sigma} \left( {\cal{H}}_\gamma \partial_z \Phi_0 - W^{\prime\prime}(\Phi_0) \Phi_1 + \partial_{zz} \Phi_1 \right) + G(\Phi_0) L_0.
\end{equation}
Multiplying by $\partial_z \Phi_0$ and integrating leads to
\begin{eqnarray}
    - \tilde{\beta} {\cal{V}} \int_{-\infty}^{+\infty} G(\Phi_0) (\partial_z \Phi_0)^2 \; dz &=&  \tilde{\sigma} \int_{-\infty}^{+\infty}
    {\cal{H}}_\gamma (\partial_z \Phi_0)^2 - \partial_z W^\prime(\Phi_0) \Phi_1 + \partial_{zz} \Phi_1 \partial_z \Phi_0 \; dz \nonumber \\
    && + L_0 \int_{-\infty}^{+\infty} G(\Phi_0) \partial_z \Phi_0 \; dz.
\end{eqnarray}
It follows $\int_{-\infty}^{+\infty} - \partial_z W^\prime(\Phi_0) \Phi_1 + \partial_{zz} \Phi_1 \partial_z \Phi_0 \; dz = \int_{-\infty}^{+\infty} \Phi_1 \partial_z (-W^\prime(\Phi_0) + \partial_{zz} \Phi_0) \; dz = 0$ and thus
\begin{equation}
    - \tilde{\beta} {\cal{V}} \int_{-\infty}^{+\infty} \!\!G(\Phi_0) (\partial_z \Phi_0)^2 \; dz =  \tilde{\sigma} {\cal{H}}_\gamma \int_{-\infty}^{+\infty}
     \!\!(\partial_z \Phi_0)^2 \; dz + L_0 \int_{-\infty}^{+\infty} \!\!G(\Phi_0) \partial_z \Phi_0 \; dz.
\end{equation}
With eq. \eqref{eq:tanh} we obtain $\partial_z \Phi_0 = \frac{1}{\sqrt{2}} (1 - \Phi_0^2)$ and thus $\int_{-\infty}^{+\infty} G(\Phi_0) (\partial_z \Phi_0)^2 \; dz = \frac{4 \sqrt{2}}{5}$, $\int_{-\infty}^{+\infty} (\partial_z \Phi_0)^2 \; dz = \frac{2 \sqrt{2}}{3}$ and $\int_{-\infty}^{+\infty} G(\Phi_0) \partial_z \Phi_0 \; dz = 2$ and therefore
\begin{equation} \label{eq:int}
    \beta {\cal{V}} = - \sigma {\cal{H}}_\gamma - 2 L_0.
\end{equation}
In order to determine $L_0$ we consider the constraint $\frac{1}{|{\cal{S}}|} \int_{\cal{S}} \phi_\epsilon \; d {\cal{S}} = \alpha$. Using eq. \eqref{eq:time} in $O(\epsilon^{-1})$ gives $0 = \int_\gamma {\cal{V}} \partial_z \Phi_0 \; d \gamma$ and as $\partial_z \Phi_0$ is independent of $s$ also $\int_\gamma {\cal{V}} \; d\gamma = 0$. Integrating eq. \eqref{eq:int} we thus obtain
\begin{equation}
    2 L_0 = - \sigma \frac{1}{|\gamma|} \int_\gamma {\cal{H}}_\gamma \, d \gamma
\end{equation}
which leads to the desired eq. \eqref{eq:MCF}. This analysis is not affected by considering $G_\eta$ instead of $G$.

\section{Graph formulations}

Before we numerically solve both models, the degenerate area preserving surface Allen-Cahn equation
\eqref{eq:AC} and \eqref{eq:AC_lambda} and the area preserving geodesic curvature flow \eqref{eq:MCF}, we reformulate them in a graph formulation.
We represent the closed curve $\gamma(t)$ on ${\cal{S}}$ as a graph of a function $h: \RR^2 \to \RR$, such that
\begin{equation}
    \gamma(t) = \{(\mathbf{X}, h(\mathbf{X}))^T : \mathbf{X} \in g(t) \}
\end{equation}
where $g(t)$ is a planar curve in $\RR^2$ and $\mathbf{X} = \mathbf{X}(l,t)$ its 1-periodic parameter form such that $g(t) = \{\mathbf{X}(l,t) = (X_1(l,t), X_2(l,t))^T, l \in [0,1] \}$. In \cite{Kolaretal_DCDSB_2017} eq. \eqref{eq:MCF} is analysed by means of the flow of $g(t)$:
\begin{equation}
    \beta V = - a H_g + b + c \frac{1}{\int_g \sqrt{1 + (\nabla h \cdot \mathbf{t}_g)^2} \; dg} \int_g {\cal{H}}_\gamma \sqrt{1 + (\nabla h \cdot \mathbf{t}_g)^2} \; dg
\end{equation}
with normal velocity $V$, coefficients $a > 0$, $b$ and $c$ given as
\begin{equation*}
a = \sigma \frac{1}{1 +(\nabla h \cdot \mathbf{t}_g)^2}, \;\;
b = \sigma \frac{\mathbf{t}_g^T \nabla^2 h \mathbf{t}_g (\nabla h \cdot \mathbf{n}_g)}{(1 +(\nabla h \cdot \mathbf{t}_g)^2(1 + |\nabla h|^2)}, \;\;
c = \sigma \sqrt{\frac{1 + (\nabla h \cdot \mathbf{t}_g)^2}{1 + |\nabla h|^2}} .
\end{equation*}
The unit tangent $\mathbf{t}_g$, unit outer normal $\mathbf{n}_g$ and curvature $H_g$ of the curve $g(t)$ given as
\begin{equation}
    \mathbf{t}_g = \frac{\partial_h \mathbf{X}}{ |\partial_h \mathbf{X}|}, \quad \mathbf{n}_g = \frac{1}{|\partial_h \mathbf{X}|} (\partial_h X_2, - \partial_h X_1)^T, \quad H_g = - \frac{1}{|\partial_h \mathbf{X}|}\partial_h \left(\frac{\partial_h \mathbf{X}}{|\partial_h \mathbf{X}|} \right) \cdot \mathbf{n}_g
\end{equation}
and the geodesic curvature ${\cal{H}}_\gamma$ of $\gamma(t)$
\begin{equation}
    {\cal{H}}_\gamma = - \frac{1}{\sqrt[3]{1 + (\nabla h \cdot \mathbf{t}_g)^2}}\left(
        \sqrt{1 + |\nabla h|^2} H_g - \frac{\mathbf{t}_g^T \nabla^2 h \mathbf{t}_g}{\sqrt{1 + |\nabla h|^2}} (\nabla h \cdot \mathbf{n}_g)
    \right).
\end{equation}
For details and the numerical realization we refer to \cite{Kolaretal_DCDSB_2017,MinBeKr:10,BeKrKri:09}.

The graph formulation for eq. \eqref{eq:AC} and \eqref{eq:AC_lambda} can be formulated as
\begin{eqnarray}
    && \hspace*{1cm} \epsilon \tilde{\beta} G_\eta(\phi) \partial_t \phi = \tilde{\sigma} \left( \epsilon \nabla\cdot\left(\left( I - \frac{(\nabla h)^2}{1+|\nabla h|^2}\right) \nabla \phi\right) - \frac{1}{\epsilon} W^\prime(\phi)\right) \\
    && + \frac{\tilde{\sigma} G_\eta(\phi)}{\int_{\Omega}\sqrt{1+|\nabla h|^2}\,d\Omega} \int_{\Omega}\frac{\sqrt{1+|\nabla h|^2}}{G_\eta(\phi)} \left(\frac{1}{\epsilon}W'(\phi)-\epsilon \nabla\cdot\left(\left( I - \frac{(\nabla h)^2}{1+|\nabla h|^2}\right) \nabla \phi\right)\right)\,d\Omega, \nonumber
\end{eqnarray}
where $\Omega \subset \RR^2$ and $h(\Omega) = {\cal{S}}$. For the numerical realization we follow \cite{DLW_JCP_2012} and introduce a relaxation rate for the Lagrange multiplier, which can also be interpreted as an additional penalization of the area of the form $c \left( \frac{1}{|\cal{S}|} \int_{\cal{S}} \phi \; d{\cal{S}} - \alpha \right)$, with penalization parameter $c$. The resulting semi-discrete graph formulation with $\phi^n = \phi(t^n)$. reads
\begin{eqnarray*}
    && \epsilon \tilde{\beta} G_\eta(\phi^n) \frac{\phi^{n+1} - \phi^n}{\tau} = \tilde{\sigma} \left( \epsilon \nabla\cdot\left(\left( I - \frac{(\nabla h)^2}{1+|\nabla h|^2}\right) \nabla \phi^{n+1}\right) - \frac{1}{\epsilon} W^\prime(\phi^{n+1})\right) \\
    && + \frac{\tilde{\sigma} G_\eta(\phi^n)}{\int_{\Omega}\sqrt{1+|\nabla h|^2}\,d\Omega} \int_{\Omega}\frac{\sqrt{1+|\nabla h|^2}}{G_\eta(\phi^n)} \left(\frac{1}{\epsilon}W'(\phi^n)-\epsilon \nabla\cdot\left(\left( I - \frac{(\nabla h)^2}{1+|\nabla h|^2}\right) \nabla \phi^{n}\right)\right)\,d\Omega \\
    && - c \left(\frac{\int_{\Omega}\sqrt{1+|\nabla h|^2} \phi^n \,d\Omega}{\int_{\Omega}\sqrt{1+|\nabla h|^2}\,d\Omega} - \alpha \right),
\end{eqnarray*}
with $W^\prime(\phi^{n+1}) \approx W^\prime(\phi^n) +
W^{\prime\prime}(\phi^n)(\phi^{n+1}-\phi^{n})$. The resulting equation for $\phi^{n+1}$ is linear and discretized in space by standard $P^1$ finite elements. The problem is implemented in AMDiS \cite{Veyetal_CVS_2007,Witkowskietal_ACM_2015} and the linear system is solved with the direct solver of UMFPACK.

\section{Numerical results}

We consider the 4 examples provided in \cite{Kolaretal_DCDSB_2017} and take the numerical solutions provided in \cite{Kolaretal_DCDSB_2017} (with $M=200$ finite volumes) as benchmark solutions. The problem settings are provided in Table \ref{tab:problem}.

\begin{table}[tbhp]
\footnotesize
\caption{Numerical examples, initial parametrization of curve and height profile. $r(l) = 1 + 0.65 \cos{(10 \pi l})$.}
\label{tab:problem}
\begin{center}
    \begin{tabular}{|c|c|c|}\hline
         & $\mathbf{X}_0$, $l \in [0,1]$ & $h$ \\ \hline
    Problem 1    & $\mathbf{X}_0 = (\frac{1}{4} + r(l) \cos{2 \pi l}, - \frac{1}{4} + r(l) \sin{2\pi l})^T$ & $h(x,y) = \sqrt{4 - x^2 - y^2}$ \\
    Problem 2    & $\mathbf{X}_0 = (\cos{2\pi l}, \frac{1}{10} + \sin{2 \pi l})^T$ & $h(x,y) = y^2$ \\
    Problem 3    & $\mathbf{X}_0 = (\cos{2\pi l}, -\frac{1}{5} + \sin{2 \pi l})^T$ & $h(x,y) = \sin{\pi y}$ \\
    Problem 4    & $\mathbf{X}_0 = (\frac{1}{2} \cos{2 \pi l}, \sin{2 \pi l})^T$ & $h(x,y) = x^2 - y^4$\\
    \hline
    \end{tabular}
\end{center}
\end{table}

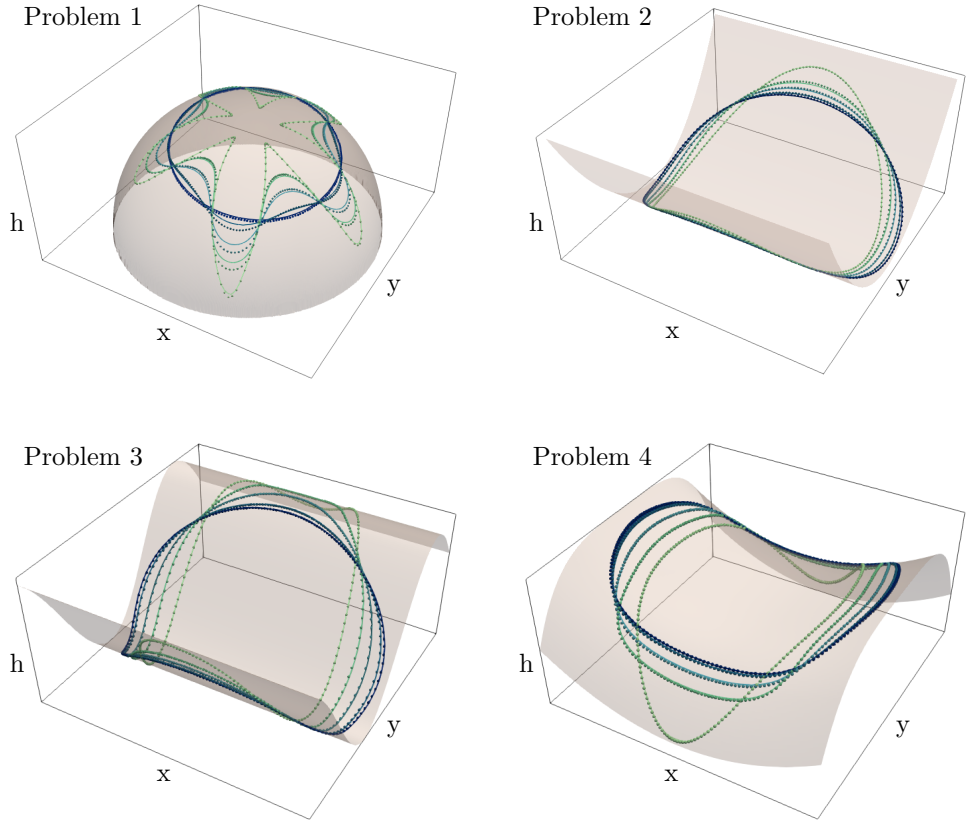
\begin{figure}
    \begin{tikzpicture}[x=0.52\textwidth, y=0.45\textwidth]
    \node[] (ex1) at (0,0)
        {\includegraphics[width=.45\textwidth]{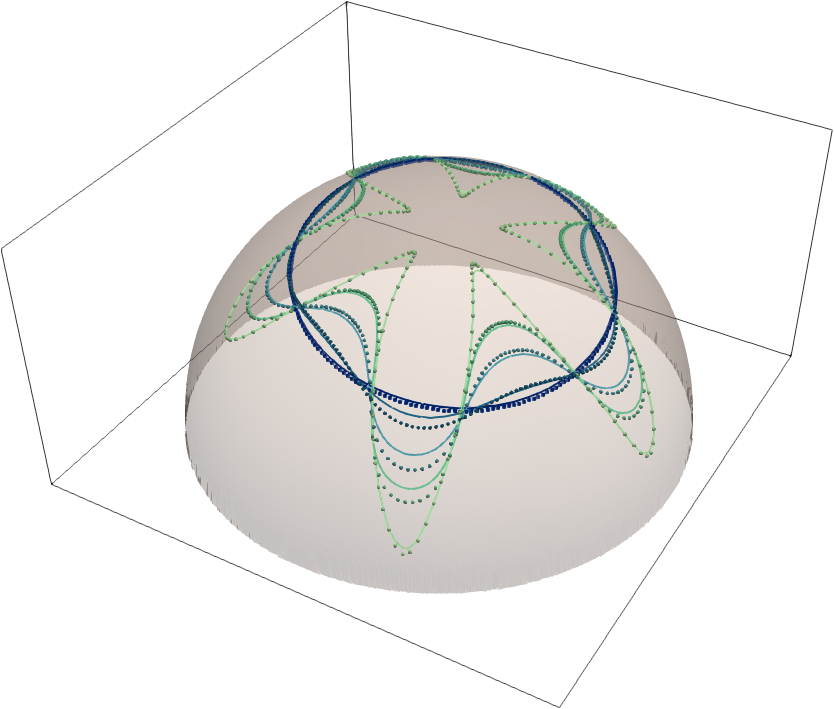}};
    \node[] at (-0.3,0.4){Problem 1};
    \node[] at (-0.43,-0.07) {h};
    \node[] at (-0.142,-0.325) {x};
    \node[] at (0.31,-0.22) {y};
    \node[] (ex2) at (1,0)
        {\includegraphics[width=.45\textwidth]{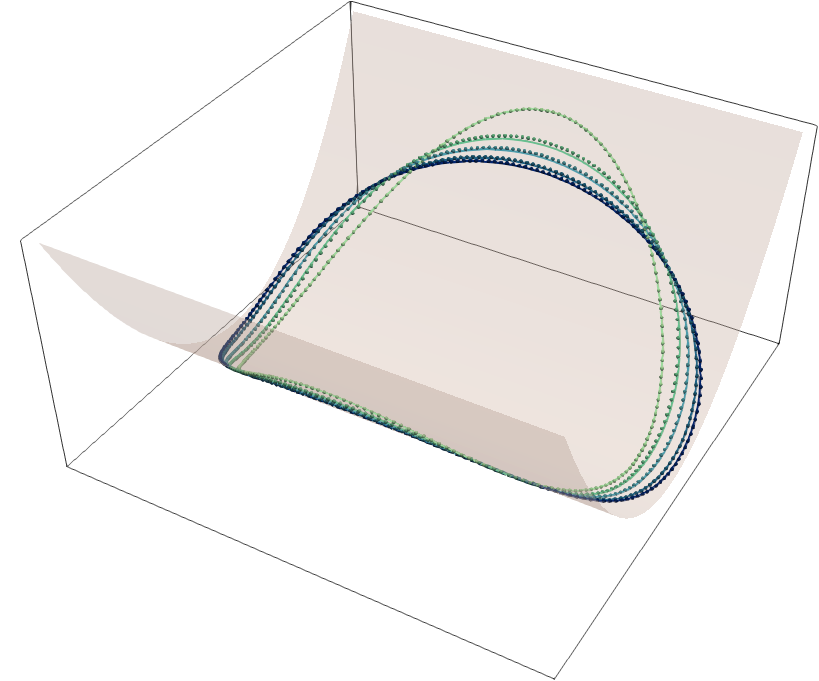}};
    \node[] at (0.7,0.4){Problem 2};
    \node[] at (0.59,-0.07) {h};
    \node[] at (1-0.143,-0.32) {x};
    \node[] at (1+0.31,-0.22) {y};
    \node[] (ex3) at (0,-1)
        {\includegraphics[width=.45\textwidth]{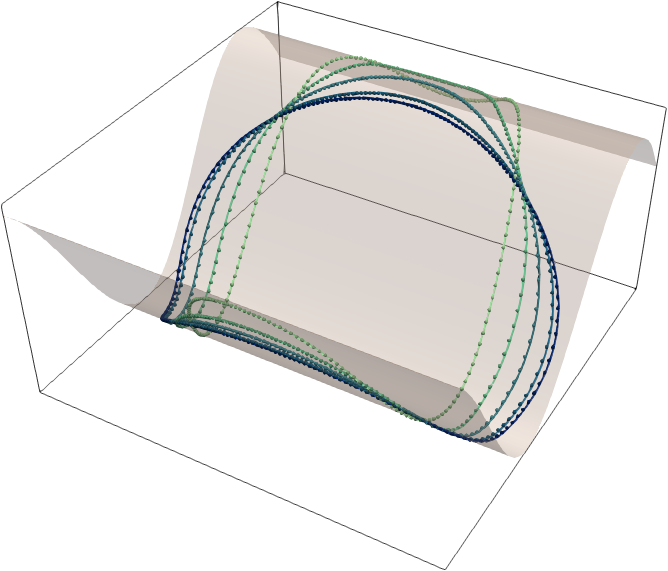}};
    \node[] at (-0.3, -0.6){Problem 3};
    \node[] at (-0.43,-1.07) {h};
    \node[] at (-0.142,-1.325) {x};
    \node[] at (0.31,-1.22) {y};
    \node[] (ex4) at (1,-1)
        {\includegraphics[width=.45\textwidth]{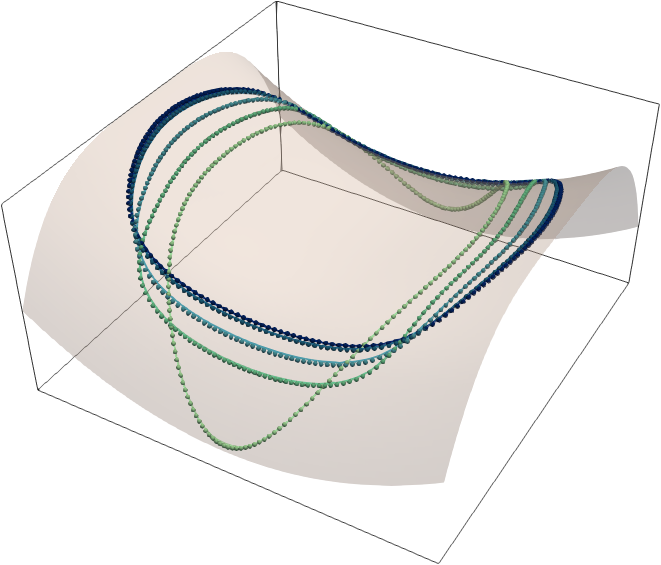}};
    \node[] at (0.7,-0.6){Problem 4};
    \node[] at (0.57,-1.07) {h};
    \node[] at (0.855,-1.325) {x};
    \node[] at (1.31,-1.22) {y};
    \end{tikzpicture}
    \caption{Comparison of the numerical solutions of the degenerate area preserving surface Allen-Cahn equation (solid) with the numerical solutions of the corresponding sharp interface limit as reference solutions (dotted). Shown is the zero contour of the phase field function $\phi_\epsilon$ for $\epsilon = 0.025$. The evolution in time is shown in color, running from light green to dark blue. The problem numbers correspond to Table \ref{tab:problem} and the time instances shown are Problem 1: $t = 0, 0.05, 0.1, 0.2, 0.4 ,1$, Problem 2: $t = 0, 0.25, 0.5, 1, 2 $, Problem 3: $t = 0, 0.5, 1, 2, 4$, and Problem 4: $t = 0, 0.25, 0.5, 1, 2.4$.}
    \label{fig:results}
\end{figure}

Figure \ref{fig:results} shows selected time instances of the solutions. The spatial resolution considers at least 10 mesh points across the projected interface in $\Omega$ and the time step is chosen to ensure the CFL condition with $\tau \approx k^2$, where $k$ is the corresponding mesh size within the diffuse interface. The mesh is adaptivly refined to ensure these conditions. Other numerical parameters are chosen as $\eta = 0.01$ and $c = 2000$. The physical parameters are set as $\beta = 1$ and $\sigma = 1$.

We also measure the space-time error of the Hausdorff distance for different $\epsilon$ of the projected curves $\gamma_{200}(t)$ (benchmark solutions) and $\gamma_\epsilon(t)$ ($\phi_\epsilon = 0$ level-sets) onto $\Omega$. We consider the $L^2$ norm in time of this distance. In addition we provide the Hausdorff distance for the reached equilibrium solutions, see \cref{tab:distance}. As the benchmark solution is also just a numerical approximation, we only discuss convergence qualitatively and do not consider any order of convergence. The values at least indicate a reduction of the considered errors. A more detailed and analytically supported convergence study of the numerical solutions requires to extend results of \cite{FP_NM_2003} to surfaces, which is beyond the scope of this paper.

\begin{figure}
    \includegraphics[width =\textwidth]{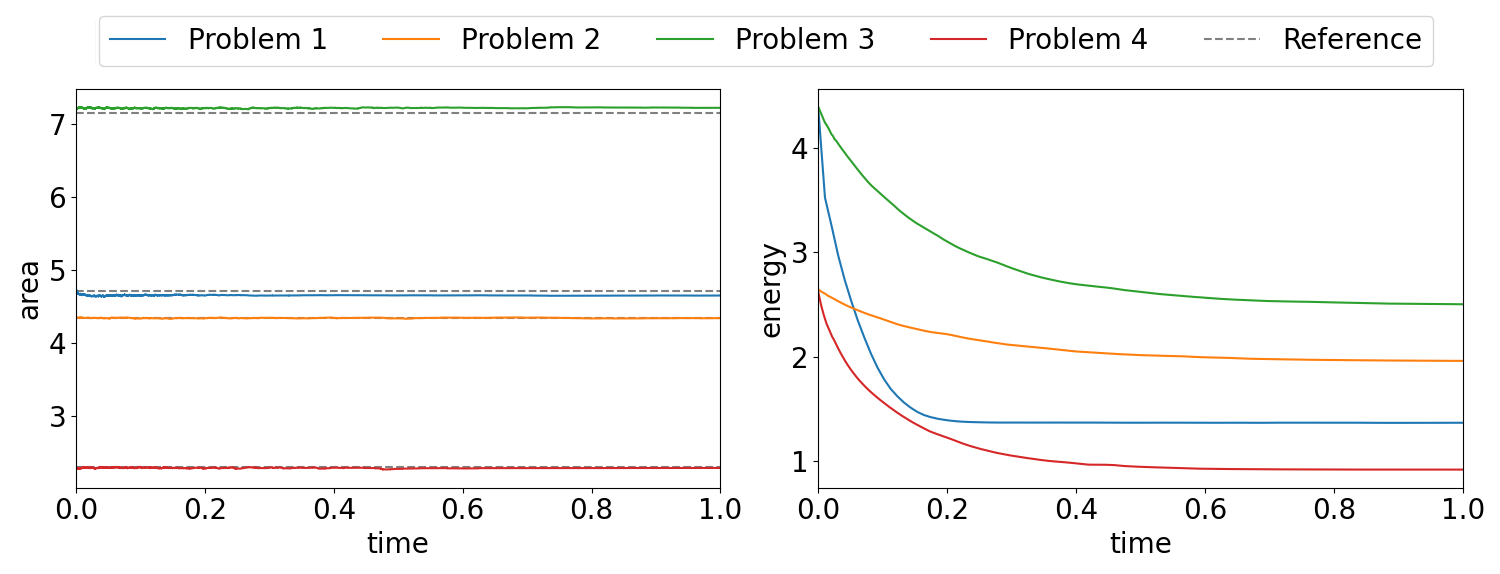}
    \caption{Time evolution of the area enclosed by the zero contour of the phase field function $\phi_\epsilon$ (left) and the system energy \cref{eq:FdGCH} (right) for $\epsilon = 0.025$. The time is normalized for all four problems. In the left plot we show the reference areas for the four problems according to \cite{Kolaretal_DCDSB_2017}. For $\epsilon = 0.025$, the maximum absolute deviation from the reference area is $0.082$ for Problem 1, $0.0180$ for Problem 2, $0.0745$ for Problem 3, and $0.0344$ for Problem 4. The mean absolute deviations are $0.0646$ for Problem 1, $0.0059$ for Problem 2, $0.0633$ for Problem 3, and $0.0133$ for Problem 4.}
    \label{fig:area_energy}
\end{figure}

Additionally we analyse the quality of area conservation and the evolution of the system energy \cref{eq:FdGCH}. The results are shown in \cref{fig:area_energy} and demonstrate the desired properties.
\begin{table}[tbhp]
\footnotesize
\caption{Hausdorff distance of the projected curves $\gamma_{200}(t)$ (benchmark solution) and $\gamma_\epsilon(t)$ ($\phi_\epsilon = 0$ level-set) onto $\Omega$ in space. Shown are the $L^2$-norm in time of the Hausdorff distances and the values at the equilibrium state for the four problems shown in \cref{tab:problem}}
\label{tab:distance}
\begin{center}
    \begin{tabular}{|c|cc|cc|cc|cc|}\hline
    &\multicolumn{2}{|c|}{\bf Problem 1}&\multicolumn{2}{|c|}{\bf Problem 2}
    &\multicolumn{2}{|c|}{\bf Problem 3}&\multicolumn{2}{|c|}{\bf Problem 4}\\\hline
    $\epsilon$     & $L^2$-norm & equil & $L^2$-norm & equil & $L^2$-norm & equil & $L^2$-norm & equil \\\hline
    0.1   & 0.0075 & 0.0577  &  0.0014 & 0.0096  &  0.0033 & 0.0211  &  0.0031 & 0.0175\\ \hline
    0.05  & 0.0043 & 0.0332  &  0.0009 & 0.0044  &  0.0019 & 0.0105  &  0.0011 & 0.0069\\ \hline
    0.025 & 0.0032 & 0.0311  &  0.0007 & 0.0026  &  0.0015 & 0.0074  &  0.0007 & 0.0056\\ \hline
    \end{tabular}
\end{center}
\end{table}

We would like to remark that the equilibrium shapes in Problem 2 and Problem 3, which are considered on ruled surfaces, are isotropic. Unrolling the surfaces provide the circular shapes of the final curves. This is consistent with known results for (reaction-)diffusion problems on curved surfaces. E.g. for the surface heat equation it is known that the heat kernel to lowest order only depends on the Gaussian curvature of the underlying surface \cite{McKean_JDE_1967,Faraudo_JCP_2002}. As this is zero for ruled surfaces, the surface should not have any effect on the evolution and a circular equilibrium shape on the surface, as in flat space, can be expected.

\section{Conclusions}

We propose a phase field approximation for area preserving geodesic curvature flow. The considered equation is a degenerate area preserving surface Allen-Cahn equation. The connection between both models is established by formal matched asymptotic analysis and confirmed by numerical solutions for different problems in a graph formulation.

 The degeneracy in the surface Allen-Cahn equation results from the de Gennes factor $G(\phi)$ in the energy \eqref{eq:FdGCH}, see \cite{SVW_MMAS_2021}. While the formal matched asymptotic analysis also holds for $G(\phi) = 1$ the numerical results relay on the de Gennes factor. It ensures $\phi \in [-1,1]$ much better than without it. This is a desired feature also in 2D and 3D, where the factor $G(\phi)$ is used in various applications \cite{Naffouti2017,AEVT_PRM_2020}. However, on curved surfaces it is even more essential as deviations have a more dramatic effect due to the spatially varying geometric properties of the surface which can enhance the resulting errors. These geometric properties also need to be considered in the mesh resolution in the graph formulation. We need to ensure a desired resolution of the projected diffuse interface. Also the additional penalization of the area \cite{DLW_JCP_2012} helps to obtain the shown convergence results. If all these aspects are considered the proposed phase field approximation provides an appropriate way to solve the highly non-linear problem of area preserving geodesic curvature flow by standard tools for solving partial differential equations in 2D.

\bibliographystyle{elsarticle-num}
\bibliography{literatur}{}

\end{document}